\documentclass[12pt]{l4dc2021}
\title{Minimax Adaptive Control for a Finite Set of Linear Systems}
\usepackage{times}
\usepackage{latexsym}
\newtheorem{prp}{Proposition}
\newtheorem{thm}[prp]{Theorem}

\newtheorem{lem}[prp]{Lemma}
\newenvironment{pf}{\smallbreak\noindent{\it Proof. }}{\hfill$\Box$\smallbreak}
\newenvironment{pf*}[1]{\smallbreak\noindent{\it #1}}{\hfill$\Box$\smallbreak}
\newcounter{definition}

\newcounter{remark}
\newenvironment{rmk}{\addtocounter{remark}{1}\smallbreak\noindent
  {\em Remark \theremark.}}{\smallbreak}
\newcommand{\realR}{\mathbb{R}}
\DeclareMathOperator{\trace}{trace}
\author{%
 \Name{Anders Rantzer} \Email{rantzer@control.lth.se}\\
 \addr Automatic Control LTH\\
 Lund University\\ 
 Box 118\\ 
 SE-221 00 Lund\\
 Sweden
}

\begin{document}

\maketitle

\begin{abstract}%
An adaptive controller with bounded $l_2$-gain from disturbances to errors 
is derived for linear time-invariant systems with uncertain parameters restricted to a finite set. The gain bound refers to the closed loop system, including the non-linear learning procedure. As a result, robustness to unmodelled dynamics (possibly nonlinear and infinite-dimensional) follows from the small gain theorem. The approach is based on a new zero-sum dynamic game formulation, which optimizes the trade-off between exploration and exploitation. An explicit upper bound on the optimal value function is stated in terms of semi-definite programming and a corresponding simple formula for an adaptive controller achieving the upper bound is given. Once the uncertain parameters have been sufficiently estimated, the controller behaves like standard $H_\infty$ optimal control. 
\end{abstract}

\begin{keywords}%
  adaptive control, real-time learning
\end{keywords}

\section{Introduction}
\begin{figure}
  \begin{center}
  \setlength{\unitlength}{.009mm}%
  \begin{picture}(9366,7266)(1318,-8494)
  \put(1001,-4500){errors}%
  \put(4000,-8461){\framebox(4200,1700){\begin{tabular}{cc}Minimax\\Adaptive Controller
      \end{tabular}}}
  \put(4000,-5611){\framebox(4200,1700){\begin{tabular}{cc}
      Linear dynamics with\\ uncertain parameters
    \end{tabular}}}
  \put(4000,-4561){\line(-1, 0){1150}}
  \put(2851,-4561){\line( 0, 1){2550}}
  \put(2851,-2011){\vector( 1, 0){1650}}
  \put(4000,-5161){\line(-1, 0){1150}}
  \put(2851,-5161){\line( 0,-1){2550}}
  \put(2851,-7711){\vector( 1, 0){1150}}
  \put(4000,-4861){\vector(-1, 0){2900}}
  \put(10651,-4861){\vector(-1, 0){2450}}
  \put(7501,-2011){\line( 1, 0){1650}}
  \put(9151,-2011){\line( 0,-1){2550}}
  \put(9151,-4561){\vector(-1, 0){950}}
  \put(8201,-7711){\line( 1, 0){950}}
  \put(9151,-7711){\line( 0, 1){2550}}
  \put(9151,-5161){\vector(-1, 0){950}}
  \put(9901,-4500){disturbances}%
  \put(4501,-2761){\framebox(3000,1500){\begin{tabular}{cc}
        Unmodelled\\ dynamics
      \end{tabular}}}
  \end{picture}%
  \end{center}
  \caption{Given a finite set of stabilizable linear models, the objective of this paper is to construct an adaptive controller that minimizes the $l_2$-gain from disturbances to errors under worst case values of the uncertain parameters. The gain bound gives a guaranteed robustness to unmodelled (possibly nonlinear and infinite-dimensional) dynamics. Optimal trade-off between exploration and exploitation is needed for minimization of the gain.}
\label{fig:blockdiag}
\end{figure}
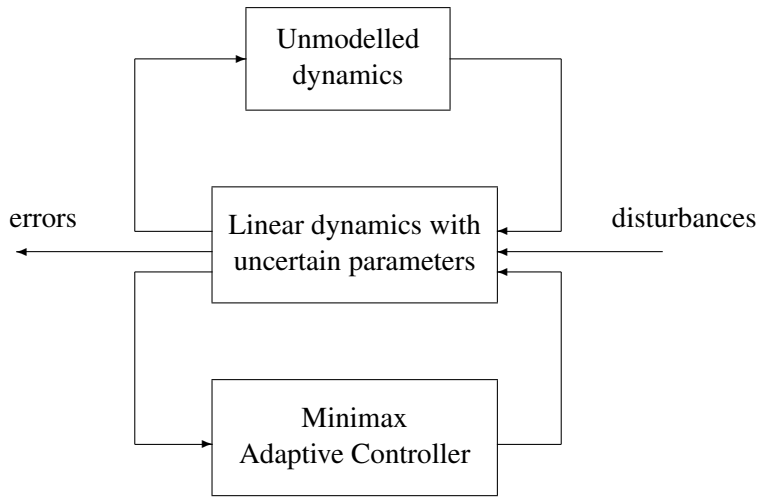

The history of adaptive control dates back at least to aircraft autopilot development in the 1950s. Following the landmark paper \cite{ast+wit73}, a surge of research activity during the 1970s derived conditions for convergence, stability, robustness and performance under various assumptions. For example, \cite{ljung77ana} analysed adaptive algorithms using averaging, \cite{goodwin1981discrete} derived an algorithm that gives mean square stability with probability one, while \cite{guo1995convergence} gave conditions for the optimal asymptotic rate of convergence. On the other hand, conditions that may cause instability were studied in \cite{ega79book}, \cite{ioannou1984instability} and \cite{rohrs1985robustness}. Altogether, the subject has a rich history documented in numerous textbooks, such as \cite{aastrom2013adaptive}, \cite{goodwin2014adaptive}, \cite{narendra2012stable}, \cite{sastry2011adaptive} and \cite{astolfi2007nonlinear}. In this paper, the focus is on worst-case models for disturbances and uncertain parameters, as discussed in \cite{cusumano1988nonlinear,sun1987theory,vinnicombe2004examples} and \cite{2004megretskinonlinear}. The ``minimax adaptive'' paradigm illustrated in Figure~\ref{fig:blockdiag} was introduced for linear systems in \cite{didinsky1994minimax} and nonlinear systems in \cite{pan1998adaptive}. 

More recently, there has been an explosion of research interest in the boundary between machine learning, system identification and adaptive control. For a review, see for example \cite{matni2019self}. Most of the studies are carried out in a stochastic setting, but recently works connecting to $H_\infty$ control have started to appear \cite{NEURIPS2018_0ae3f79a}. There is also recent work connecting adversarial learning and minimax regret analysis to robust control \cite{zhang2020stability}. 
However, it is worth noting that unlike most recent contributions the minimax approach of this paper will not assume that a stabilizing controller is known in advance.

This paper is an extension and generalization of the approach presented in \cite{rantzer2020ifac}. The outline is as follows: After some notation in sections~2, we give the problem formulation in section~3, followed by preliminary results on minimax dynamic programming in section~4.
The main result is first presented in section~5 for the special case of input sign uncertainty with a double integrator as an illustrating example, then proved for the general case in section~6. Some calculations are deferred to an appendix.

\section{Notation}
The set of $n\times m$ matrices with real coefficients is denoted $\realR^{n\times m}$. The transpose of a matrix $A$ is denoted $A^\top$. For a symmetric matrix $A\in\realR^{n\times n}$, we write $A\succ0$ to say that $A$ is positive definite, while $A\succeq0$ means positive semi-definite. For $A, B\in\realR^{n\times m}$, the expression $\langle A,B\rangle$ denotes the trace of $A^\top B$. Given $x\in\realR^n$ and $A\in\realR^{n\times n}$, the notation $|x|^2_A$ means 
$x^\top Ax$. Similarly, given $B\in\realR^{m\times n}$ and $A\in\realR^{n\times n}$, the trace of $B^\top AB$ is denoted $\|B\|^2_A$.

\section{Minimax Adaptive Control}
\label{sec:minmax}

Let $Q\in\realR^{n\times n}$ and $R\in\realR^{m\times m}$ be positive definite. Given a compact set $\mathcal{M}\subset\realR^{n\times n}\times\realR^{n\times m}$ and a number $\gamma>0$, we are interested to compute 
\begin{align}
  J_*(x_0)&:=\inf_\mu\underbrace{\sup_{w,A,B,N}\sum_{t=0}^N\left(|x_t|_Q^2+|u_t|_R^2-\gamma^2|w_t|^2\right)}_{J_\mu(x_0)},
\label{eqn:infsup}
\end{align}
where $(A,B)\in\mathcal{M}$, $w_t\in\realR^n$, $N\ge0$ and the sequences $x$ and $u$ are generated according to
\begin{align}
  x_{t+1}&=Ax_t+Bu_t+w_t& t&\ge0 
\label{eqn:plant}\\
  u_t&=\mu_t(x_0,\ldots,x_t,u_0,\dots,u_{t-1}).
\label{eqn:mu_LQ}
\end{align}
The problem above can be viewed as a dynamic game, where the $\mu$-player tries to minimize the cost, while the $(w,A,B)$-player tries to maximize it. If it wasn't for the uncertainty in $(A,B)$, this would be the standard game formulation of $H_\infty$ optimal control \cite{Basar/B95}. In our formulation, the maximizing player can choose not only $w$, but also the matrices $A,B$. The matrices are unknown but constant, so an optimal feedback law tends to ``learn'' $A$ and $B$ early on, in order to exploit this knowledge later. 
Such nonlinear adaptive controllers can stabilize and optimize the behavior also when no linear time-invariant controller can simultaneously stabilize (\ref{eqn:plant}) for all $(A,B)\in\mathcal{M}$. 

To accommodate the uncertainty in $(A,B)$ when deciding $u_t$, it is natural for the controller to consider historical data collected in the matrix
\begin{align}
  Z_t&=\sum_{\tau=0}^{t-1}{\small\begin{bmatrix}-x_{\tau+1}\\x_\tau\\u_\tau\end{bmatrix}
  \begin{bmatrix}-x_{\tau+1}\\x_\tau\\u_\tau\end{bmatrix}^\top}
\label{eqn:Z_t}
\end{align}
since this gives
$
  \left\|\begin{bmatrix}I\;\;\,A\;\;\,B\end{bmatrix}^\top\right\|^2_{Z_t}
  =\sum_{\tau=0}^{t-1}|w_\tau|^2
$. In fact, we will prove in Theorem~\ref{thm:intro} that knowledge of $x_t$ and $Z_t$ is sufficient for optimal control. The result will be based on the following reformulated problem:
Given $Q\succ0, R\succ0$, $\gamma>0$, a set $\mathcal{M}$ and the system
\begin{align}
  \left\{\begin{array}{ll}
    x_{t+1}=v_t\\
    Z_{t+1}=Z_t+{\small\begin{bmatrix}-v_t\\x_t\\u_t\end{bmatrix}
             \begin{bmatrix}-v_t\\x_t\\u_t\end{bmatrix}^\top},\qquad &Z_0=0,
  \end{array}\right.
\label{eqn:xZ}
\end{align}
find a control law 
\begin{align}
  u_t&=\eta(x_t,Z_t)
\label{eqn:eta}
\end{align}
to minimize
{\begin{align}
  \inf_\eta\sup_{v,N}\Bigg[\sum_{t=0}^N\left(|x_t|_Q^2+|u_t|_R^2\right)
  -\gamma^2\!\!\!\min_{(A,B)\in\mathcal{M}}\left\|\begin{bmatrix}I&A&B\end{bmatrix}^\top\right\|^2_{Z_{N+1}}\Bigg]
\label{eqn:infsup2}
\end{align}
}when $x,u,Z$ are generated from $v$ 
according to
 (\ref{eqn:xZ})-(\ref{eqn:eta}).
\smallskip

In this formulation, the unknown $(A,B)$ does not appear in the dynamics, only in the penalty of the final state. As a consequence, no past states are needed in the control law (\ref{eqn:eta}), only the current state $(x_t,Z_t)$. In fact, the problem is a standard zero-sum dynamic game \cite{basar1999dynamic}, which will next be addressed by dynamic programming.

\section{Minimax Dynamic Programming}
Define the operators $\mathcal{F}$ and $\mathcal{F}_u$ by
\begin{align*}
  {\mathcal{F}}V(x,Z):=
     &\min_u\underbrace{\max_v\left\{|x|_Q^2+|u|_R^2
     +V\left(v,Z+{\small\begin{bmatrix}-v\\x\\u\end{bmatrix}
     \begin{bmatrix}-v\\x\\u\end{bmatrix}^\top}\right)\right\}}_{\mathcal{F}_uV(x,Z)}.
\end{align*}
Then the following result holds:

\begin{thm}
  Given $Q,R,\mathcal{M}$, define the operator ${\mathcal{F}}$ as above and $V_0,V_1,V_2\ldots$ according to the iteration\\[-2mm]
  \begin{align}
    V_0(x,Z)&=-\gamma^2\min_{(A,B)\in\mathcal{M}}\big\|\begin{bmatrix}I\;\;A\;\;\,B\end{bmatrix}^\top\big\|^2_Z\label{eqn:Bellman_init}\\
    V_{k+1}(x,Z)&={\mathcal{F}}V_k(x,Z).
  \label{eqn:Bellman_t}
  \end{align}
  The expressions (\ref{eqn:infsup}) and (\ref{eqn:infsup2}) have the same value. The value is finite if and only if the sequence $\{V_k(x,0)\}_{k=0}^\infty$ is upper bounded. If so, the limit $V_*:=\lim_{k\to\infty}V_k$ exists and $J_*(x_0)=V_*(x_0,0)$. 
  Defining $\eta(x,Z)$ as the minimizing value of $u$ in the expression for ${\mathcal{F}}V_*(x,Z)$ gives an optimal $\eta^*$ for (\ref{eqn:infsup2}), while the control law $\mu^*$ defined by\\[-2mm]
  \begin{align}
    \mu^*_t(x_0,\ldots,x_t,u_0,\dots,u_{t-1})
    &:=\eta^*\left(x_t,\sum_{\tau=0}^{t-1}{\small\begin{bmatrix}-x_{\tau+1}\\x_\tau\\u_\tau\end{bmatrix}
      \begin{bmatrix}-x_{\tau+1}\\x_\tau\\u_\tau\end{bmatrix}^\top}\right)
  \label{eqn:muopt}
  \end{align}
  is optimal for (\ref{eqn:infsup}). Moreover, if there exists a function $\bar{V}$ satisfying $\bar{V}\ge V_0$ and ${\mathcal{F}}_{\bar\eta}\bar{V}\le\bar{V}$, then 
  the control law $\bar\mu$ defined by $u_t=\bar\eta(x_t,Z_t)$ satisfies $J_{\bar\mu}(x_0)\le\bar{V}(x_0,0)$.
\label{thm:intro}
\end{thm}

\begin{pf}
  \begin{align*}
    V_1(x,Z)&\ge\min_u\max_v
    V_0\left(v,Z+{\small\begin{bmatrix}-v\\x\\u\end{bmatrix}
    \begin{bmatrix}-v\\x\\u\end{bmatrix}^\top}\right)\\
    &=-\gamma^2\max_u\min_{A,B,v}\left(|Ax+Bu-v|^2
    +\big\|\begin{bmatrix}I\;\;A\;\;\,B\end{bmatrix}^\top\big\|^2_Z\right)\\
    &=V_0(x,Z)
  \end{align*}
  so $V_1\ge V_0$ and the sequence $V_0,V_1,V_2,\ldots$ is monotonically non-decreasing.

  For any fixed $N\ge0$, the value of (\ref{eqn:infsup}) is bounded below by the expression
  \begin{align}
    \inf_\mu\sup_{w,A,B}\sum_{t=0}^N\left(|x_t|_Q^2+|u_t|_R^2-\gamma^2|w_t|^2\right),
  \label{eqn:infsup_N} 
  \end{align}
  where $(A,B)\in\mathcal{M}$, $w_t\in\realR^n$ and the sequences $x$ and $u$ are generated according to (\ref{eqn:plant})-(\ref{eqn:mu_LQ}). 
  The value of (\ref{eqn:infsup_N}) grows monotonically with $N$ and (\ref{eqn:infsup}) is obtained in the limit. A change of variables with $v_t:=x_{t+1}$ and $Z_t$ given by (\ref{eqn:Z_t}) shows that (\ref{eqn:infsup_N}) is equal to
  {\begin{align}
    \inf_\mu\sup_{v}\left[\sum_{t=0}^N\left(|x_t|_Q^2+|u_t|_R^2\right)-\gamma^2\min_{A,B}\big\|\begin{bmatrix}I\;\;A\;\;\,B\end{bmatrix}^\top\big\|^2_{Z_{N+1}}\right]
  \label{eqn:sup_vT}
  \end{align}
  }where $x,Z,u$ are generated by (\ref{eqn:xZ}) combined with (\ref{eqn:mu_LQ}). 
Standard dynamic programming shows that the value of (\ref{eqn:sup_vT}) is $V_{N+1}(x_0,0)$, where $V_k$ is defined by (\ref{eqn:Bellman_init})-(\ref{eqn:Bellman_t}).
  This proves that (\ref{eqn:infsup}) has a finite value if and only if the sequence $\{V_k(x,0)\}_{k=0}^\infty$ is upper bounded. The limit $V_*(x_0,0)$ is equal to the value of (\ref{eqn:infsup}). The relationship $V_k(x,Z)\le V_k(x,0)$ for $Z\succeq0$ shows that the limit then exists for all $Z\succeq0$.
      
  If (\ref{eqn:infsup2}) is finite, then $V_k$ is bounded above by (\ref{eqn:infsup2}), so also $V_*:=\lim_{k\to\infty}V_k$ is finite. Conversely, if $V_0\le\bar{V}\le\infty$ and ${\mathcal{F}}_{\bar{\eta}}\bar{V}\le\bar{V}$, we may define the sequence $W_0,W_1,W_2,\ldots$ recursively by $W_0:=V_0$ and
  \begin{align*}
    W_{k+1}(x,Z)
    &:=\mathcal{F}_{\bar\eta(x,Z)}W_k(x,Z).
  \end{align*}
  By dynamic programming,
  \begin{align*}
    W_N(x,0)
    &=\max_{v}\left[\sum_{t=0}^N\left(|x_t|_Q^2+|u_t|_R^2\right)-\gamma^2\min_{A,B}\big\|\begin{bmatrix}I\;\;A\;\;\,B\end{bmatrix}^\top\big\|^2_{Z_{N+1}}\right]
  \end{align*}
  where $x,Z,u$ are generated by (\ref{eqn:xZ}) and (\ref{eqn:eta}). Hence (\ref{eqn:infsup2}) is bounded above by $\lim_{k\to\infty}W_k(x_0,0)$. The definitions of $V_k$ and $W_k$ give by induction $\bar{V}\ge W_k\ge V_k$ for all $k$, so $V_*\le\lim_{k\to\infty}W_k\le\bar{V}$. This proves that $V_*(x_0,0)\le J_{\bar{\mu}}(x_0)\le \bar{V}(x_0,0)$ when the control law $\bar\mu$ is defined by $u_t=\bar\eta(x_t,Z_t)$. In particular, the control law is optimal if $\bar{V}=V_*$. 
\end{pf}


We conclude the section by pointing out that even though all our results assume full state measurements, they are immediately applicable also to input-output models. For example, the input-output model
\begin{align}
  y_t=-a_1y_{t-1}-\cdots-a_ny_{t-n}+b_1u_{t-1}+\cdots+b_nu_{t-n}
\label{eqn:IO}
\end{align}
has the (non-minimal) state realization $x_{t+1}=Ax_t+Bu_t$ where
\begin{align}
  x_t&=\begin{bmatrix}
    y_{t-1}\\\vdots\\y_{t-n}\\u_{t-1}\\\vdots\\u_{t-n}
  \end{bmatrix}&
  A&=\begin{bmatrix}
    -a_1&\ldots&-a_n&b_1&\ldots&b_n\\
    1\\
    &\ddots\\
    &&1\\
    0&\ldots&0&0&\ldots&0\\
    &&&1\\
    &&&&\ddots\\
    &&&&&1
  \end{bmatrix}&
  B&=\begin{bmatrix}
    0\\0\\\vdots\\0\\1\\0\\\vdots\\0
  \end{bmatrix}
\label{eqn:IOstate}
\end{align}
where all states are known from the inputs and outputs. This also shows that it is not really restrictive to consider $B$ as known. Example~1 will use this state realization for a double integrator where only the output of the second integrator is available for measurement. 

\section{Systems with unknown input direction}

The following result gives an explicit expression for an adaptive controller satisfying a pre-specified bound on the $L_2$-gain for any stabilizable model set of the form
\begin{align*}
  \mathcal{M}&:=\{(A,B),(A,-B))\}\subset\realR^{n\times n}\times\realR^{n\times m}.
\end{align*} 

\begin{thm}
  Consider $A\in\realR^{n\times n}$, $B\in\realR^{n\times m}$ and positive definite $Q\in\realR^{n\times n}$ and $R\in\realR^{m\times m}$. Suppose 
  there exists $K\in\realR^{m\times n}$ and symmetric  $P,T\in\realR^{n\times n}$ satisfying 
  $0\prec P\prec T\prec \gamma^2I$ and
  {\begin{align}
    P&\succeq Q+K^\top RK+(A-BK)^\top (P^{-1}-\gamma^{-2}I)^{-1}(A-BK)\label{eqn:P}\\
    T&\succeq Q+K^\top RK+(A+BK)^\top (P^{-1}-\gamma^{-2}I)^{-1}(A+BK)\label{eqn:T}\\
    T&\succeq Q+K^\top(R-\gamma^2B^\top B)K
        +A^\top(T^{-1}-\gamma^{-2}I)^{-1}A.\label{eqn:K}
  \end{align}
  }Then the Bellman inequality $\bar{V}\le\mathcal{F}\bar{V}$ is satisfied by
  {\small\begin{align}
    \bar{V}(x,Z):=
    &\max\left\{|x|^2_{P}
    -\gamma^2\left\|\begin{bmatrix}I&A&\pm B\end{bmatrix}^\top\right\|_Z^2,
    |x|^2_{T}-\gamma^2\trace\left(\begin{bmatrix}
      [I\;\;A]^\top[I\;\;A]&0\\0&B^\top B
    \end{bmatrix}Z\right)\right\}
  \label{eqn:Vbar}
  \end{align}
  }and the bound $J_\mu(x_0)\le|x_0|^2_{T}$ is valid for the control law (\ref{eqn:mu_LQ}) defined by
  \begin{align}
    u_t=\begin{cases}
      -Kx&\hbox{if }\sum_{\tau=0}^{t-1}(x_{\tau+1}-Ax_\tau)^\top Bu_\tau\ge0\\
      Kx&\hbox{otherwise.}
    \end{cases}
  \label{eqn:Lcontrol}
  \end{align}
\label{thm:inputsign}
\end{thm}

\begin{pf}
  Theorem~\ref{thm:inputsign} is a special case of Theorem~\ref{thm:main}, to be proved later. Defining
  \begin{align*}
    (A_1,B_1)&=(A,B)&
    (A_2,B_2)&=(A,-B)\\
    (P_{11},K_1)&=(P,K)&    
    (P_{22},K_2)&=(P,-K)\\
    P_{12}&=T&
    P_{21}&=T
  \end{align*}
  turns (\ref{eqn:Pik})-(\ref{eqn:Kcontrol}) into (\ref{eqn:P})-(\ref{eqn:Lcontrol}).
\end{pf}

\begin{example} 
See Figure~\ref{fig:doubleblock} for an illustration. Following (\ref{eqn:IO})-(\ref{eqn:IOstate}), we model a double integrator with unknown sign of the input gain as
\begin{align*}
  x_{t+1}&=\underbrace{\begin{bmatrix}
    2&-1&1\\
    1&0&0\\
    0&0&0
  \end{bmatrix}}_{A}x_t
  \pm\underbrace{\begin{bmatrix}
    0\\0\\1
  \end{bmatrix}}_Bu_t+w_t
\end{align*}
where the state is $x_t=\begin{bmatrix}
 y_t&y_{t-1}&u_{t-1}
\end{bmatrix}^\top$. Theorem~\ref{thm:main} can be applied with $\mathcal{M}=\{(A,\pm B)\}$. By first solving the Riccati equation for $P$ and $K$, then solving the matrix inequalities for $T$ by convex optimization, we get
\begin{align*}
  P&=\left[\begin{array}{rrr}
    20.61 & -11.09 &  11.09\\
      -11.09 &   7.83 &  -6.83\\
       11.09 &  -6.83 &   7.83
  \end{array}\right]&
  T&=\left[\begin{array}{rrr}
    155.0 & -84.4 & 84.4\\
   -84.4 & 89.0 & -87.5\\
    84.4 & -87.5 & 89.0
  \end{array}\right]\\
  \gamma&=19&K&=\left[\begin{array}{rrr}
    \;\;1.786 &  -1.288 &   1.288
  \end{array}\right]
\end{align*}
Simulations with the control law (\ref{eqn:Lcontrol}) are shown in Figure~\ref{fig:double}, first for $w=0$, then for $w$ being a white noise sequence. The case with worst case disturbance $w$ is not plotted, since it actually makes the state grow exponentially according to the formula $x_{t+1}=(I-\gamma^{-2}T)^{-1}Ax_t$. This does not mean that the system is unstable. Instead, what happens is that the worst disturbance $w$ grows exponentially, to keep the current state large compared to past data and prevent the adaptive controller from becoming confident about the input matrix sign. The construction can be used to verify that the 
the $l_2$-gain from $w$ to $(Q^{1/2}x,R^{1/2}u)$
when using the adaptive control law (\ref{eqn:Lcontrol}) is somewhere between $\sqrt{\|T\|}=16.8$ and $\gamma=19$. 

\begin{figure}
  \begin{center}
    \setlength{\unitlength}{.01mm}%
    \begin{picture}(10244,4044)(2079,-4883)
    \thicklines
    \put(3301,-3061){\framebox(1800,1200){Integrator}}
    \put(4501,-4861){\framebox(3700,1200){Unmodelled dynamics}}
    \put(7301,-3061){\framebox(1800,1200){Integrator}}
    \put(11301,-3061){\framebox(1200,1200){$\pm 1$}}
    \put(3301,-2461){\vector(-1, 0){2200}}
    \put(7301,-2461){\vector(-1, 0){2200}}
    \put(11301,-2461){\vector(-1, 0){2200}}
    \put(14001,-2461){\vector(-1, 0){1500}}
    \put(2501,-2461){\line( 0,-1){1800}}
    \put(10101,-4261){\vector( 0,1){1720}}
    \put(2501,-4261){\vector( 1, 0){2000}}
    \put(8201,-4261){\line( 1, 0){1900}}
    \put(10101,-2461){\circle{212}}
    \put(13501,-2161){$u_t$}
    \put(10501,-2161){$\pm u_t$}
    \put(9501,-2161){$\nu_t$}
    \put(6101,-2161){$\xi_t$}
    \put(1501,-2161){$y_t$}
    \put(7001,-1561){$\xi_{t+1}=\xi_t+\nu_t$}
    \put(3001,-1561){$y_{t+1}=y_t+\xi_t$}
    \end{picture}%
  \end{center}
\caption{No linear feedback can stabilize a double integrator with unknown input sign. However, 
Example~1 illustrates how a minimax adaptive controller can do this and even get robustness to a certain amount of unmodelled dynamics. Only the output of the second integrator is measured, but using (\ref{eqn:IO})-(\ref{eqn:IOstate}) it is still possible to write a double integrator on state space form with a fully accessible state.}
\label{fig:doubleblock}
\end{figure}
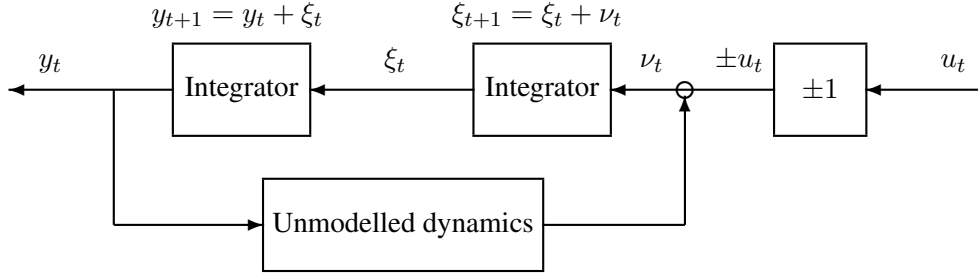

\begin{figure}[t]
\begin{center}
  \includegraphics[width=.49\hsize,height=.35\vsize]{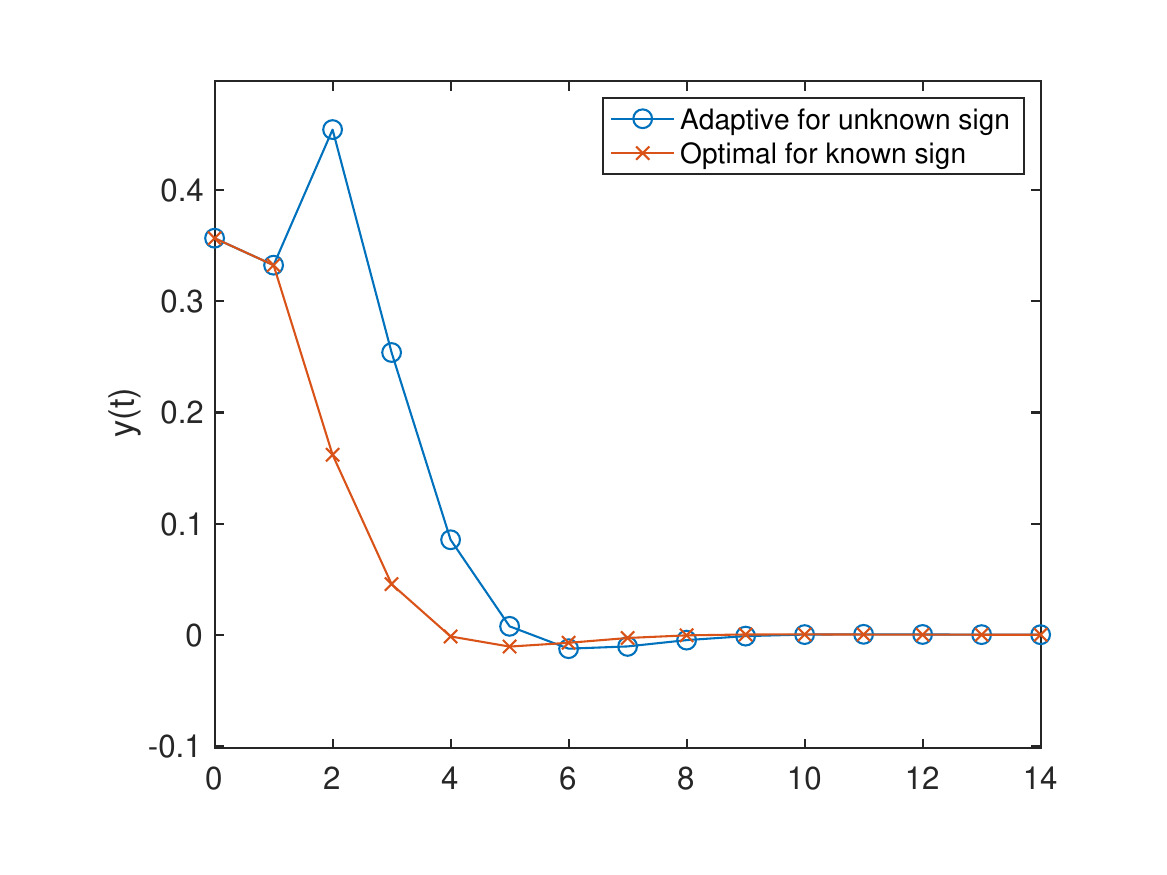}
  \includegraphics[width=.49\hsize,height=.35\vsize]{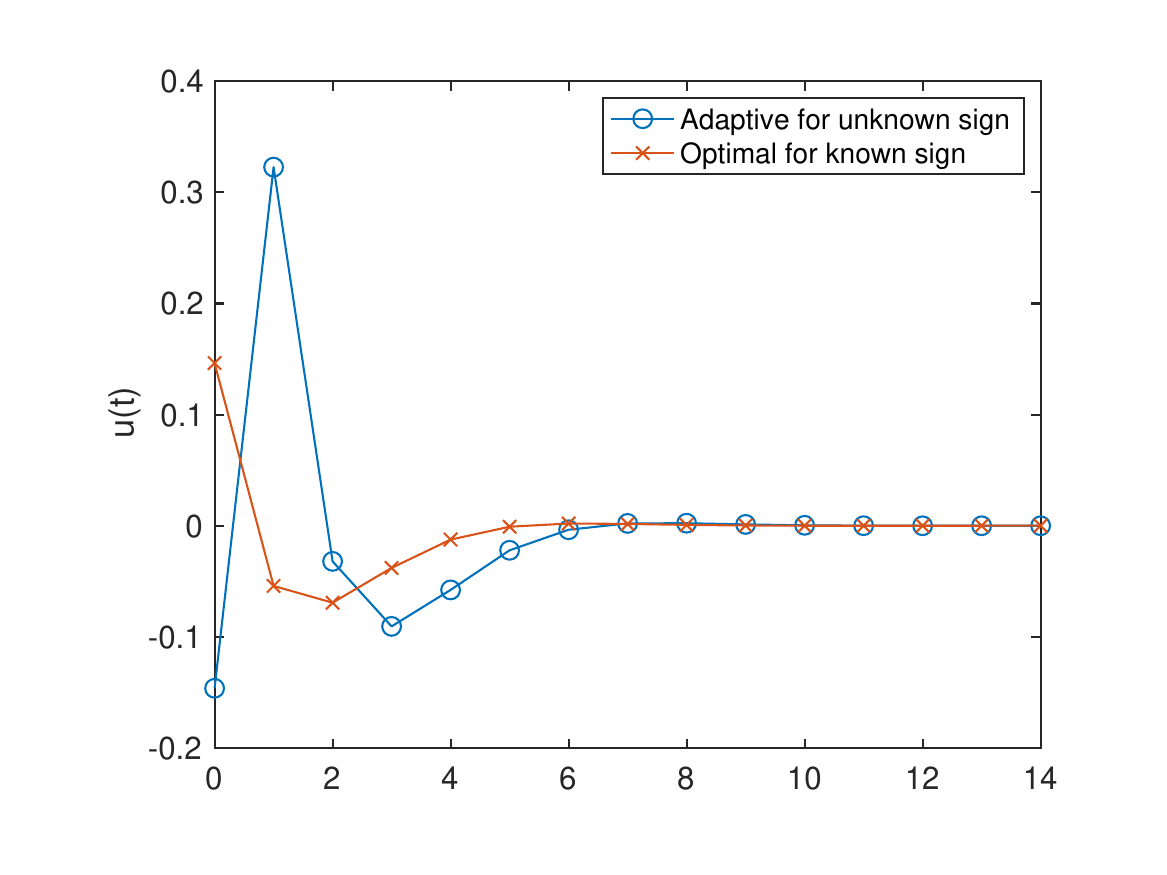}\\
  \;\;\includegraphics[width=.46\hsize,height=.35\vsize]{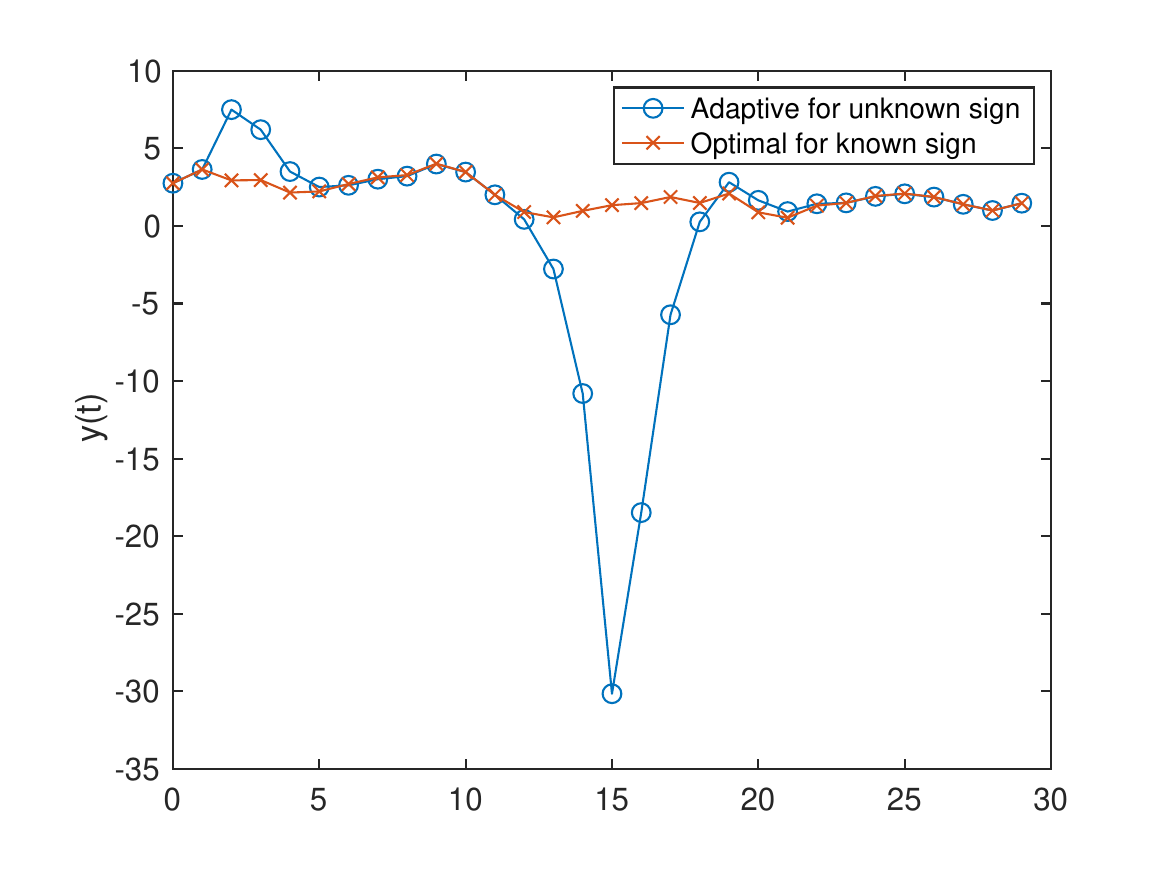}\;\quad
  \includegraphics[width=.46\hsize,height=.35\vsize]{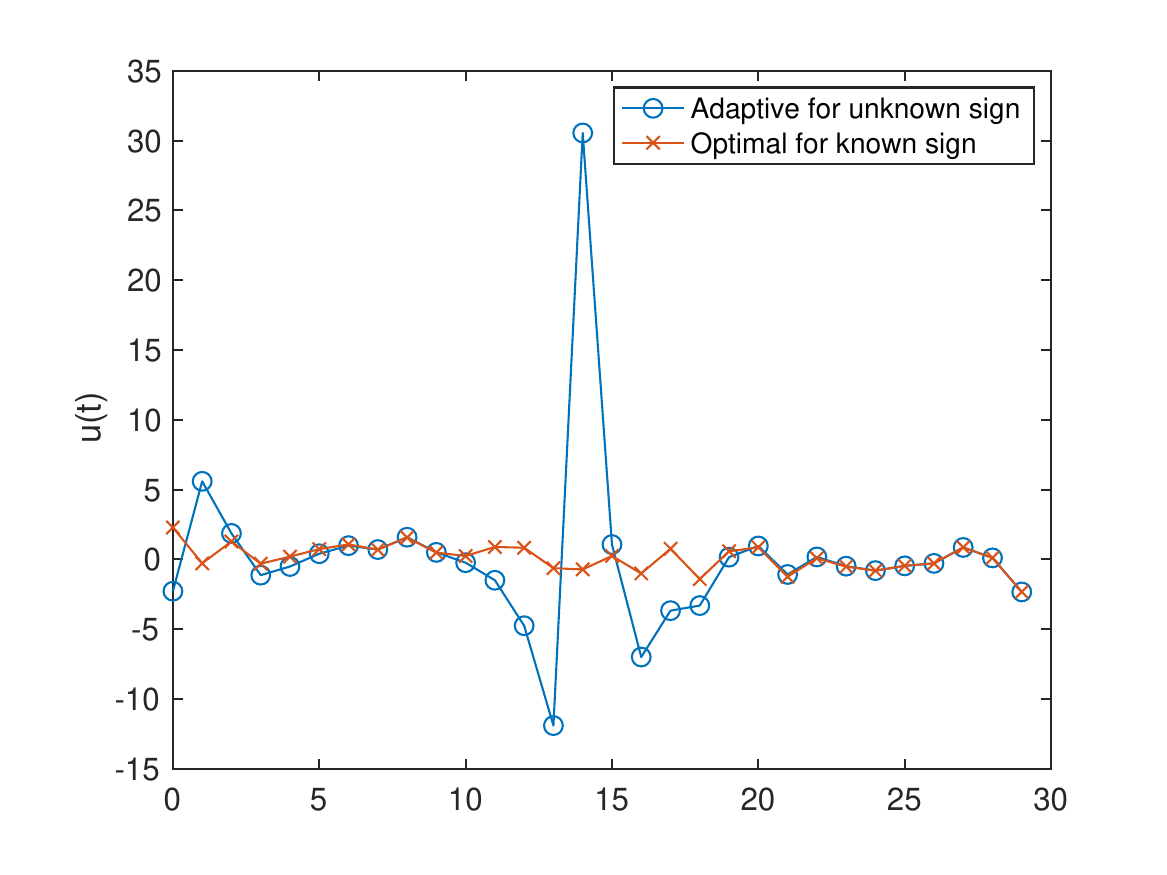}
\end{center}
\caption{Output (left) and input (right) are plotted for the double integrator with input matrix sign uncertainty, controlled using a minimax adaptive controller (blue plots with circles). For comparison, trajectories of the optimal controller with known input matrix sign are also given (red plots). The top plots have been generated with $w=0$. The minimax controller automatically increases control activity for the purpose of exploration in the beginning. In the lower plots, $w$ is a white noise sequence and the sign of the input matrix is changed at $t=10$. The adaptation takes longer the second time, since data collected before $t=10$ makes the adaptive controller reluctant to change its internal model.}
\label{fig:double}
\end{figure}
\end{example}

\section{Explicit solution of the Bellman inequality in the general case}

The following result gives an explicit expression for an adaptive controller satisfying a pre-specified bound on the $L_2$-gain for model sets of the form
\begin{align}
  \mathcal{M}&:=\{(A_1,B_1),\ldots,(A_N,B_N))\}\subset\realR^{n\times n}\times\realR^{n\times m}.
\label{eqn:Mfinite)}
\end{align} 

\begin{thm}
  Given $A_1,\ldots,A_N\in\realR^{n\times n}$, $B_1,\ldots,B_N\in\realR^{n\times m}$ and positive definite $Q\in\realR^{n\times n}$, $R\in\realR^{m\times m}$, suppose there exist $K_1,\ldots,K_N\in\realR^{m\times n}$ and $P_{ij}\in\realR^{n\times n}$ with
    $0\prec P_{ij}=P_{ji}\prec \gamma^2I$ and
  {\begin{align}
|x|^2_{P_{jk}}\ge|x|^2_Q+|K_{k}x|^2_R&+\bigg|(A_i-B_iK_k+A_j-B_jK_k)x/2\bigg|^2_{(P_{ij}^{-1}-\gamma^{-2}I)^{-1}}\notag\\
&-\gamma^2|(A_i-B_iK_k-A_j+B_jK_k)x/2|^2
\label{eqn:Pik}
\end{align}
  }for $x\in\realR^n$ and $i,j,k\in\{1,\ldots,N\}$ except if $i\ne j=k$. Then the Bellman inequality $\bar{V}\le\mathcal{F}\bar{V}$ holds for
    {\begin{align}
      \bar{V}(x,Z):=
            &\max_{i,j}\left\{|x|^2_{P_{ij}}
            -\frac{\gamma^2}{2}\left\|\begin{bmatrix}I&A_i& B_i\end{bmatrix}^\top\right\|_Z^2
                        -\frac{\gamma^2}{2}\left\|\begin{bmatrix}I&A_j& B_j\end{bmatrix}^\top\right\|_Z^2
                                    \right\}
          \label{eqn:Vbar_is}
    \end{align}
    }and the bound $J_\mu(x_0)\le\max_{i,j}|x_0|^2_{P_{ij}}$ is valid for the control law $\mu$ defined by
  \begin{align}
    u_t&=-K_{k_t}x_t,\hbox{ where }
    k_t=\arg\min_{i}\sum_{\tau=0}^{t-1}|A_ix_\tau+B_iu_\tau-x_{\tau+1}|^2.\label{eqn:Kcontrol}
  \end{align}
\label{thm:main}
\end{thm}

\begin{rmk}
  One way to approach the solution of (\ref{eqn:Pik}) is to first solve the Riccati equations
  \begin{align*}
    |x|^2_{P_{ii}}&=\min_u\max_w\left[|x|^2_Q+|u|^2_R-\gamma^2|w|^2+|A_ix+B_iu+w|^2_{P_{ii}}\right]
  \end{align*}
  for ${P_{ii}}$ and the minimizing control laws $u=-K_ix$, then use semi-definite programming to determine $P_{ij}$ for $i\ne j$. If no solution exists, then increase $\gamma$. However, it should be noted that this approach is not guaranteed to find a solution with the smallest possible $\gamma$.
\end{rmk}

\goodbreak

\begin{rmk}
  The inequalities (\ref{eqn:Pik}) have interesting interpretations. First, for $i=j=k$ it is the standard $H_\infty$ Riccati inequality, which specifies that $|x|^2_{P_{ii}}$ is an upper bound on the minimax cost function for the system $(A_i,B_i)$ when the model is known. Second, for $i=j\ne k$ it
  shows that the difference in step cost due to use of a controller that does not match the model is upper bounded by $|x|^2_{P_{ik}}-|x|^2_{P_{ii}}$. Third, (\ref{eqn:Pik}) also shows that if improved model knowledge is taken into account, then cost-to-go is decreasing even when the ``wrong'' controller is used. This is due to the last term of the inequality, which represents improved knowledge about the system due to learning. 
\end{rmk}

\begin{pf}
  Define for $i,j\in\{1,\ldots,N\}$
  \begin{align*}
    V^{ij}(x,Z)&:=|x|^2_{P_{ij}}-(z_i+z_j)/2&
    z_i&:=\gamma^2\left\|\begin{bmatrix}I&A_i&B_i\end{bmatrix}^\top\right\|_Z^2
  \end{align*}
  Then $\bar{V}=\max_{i,j}V^{ij}$. Define $k:=\arg\min_iz_i$. If $i=j=k$, then
  \begin{align*}
    \mathcal{F}_{-K_kx}V^{ii}(x,Z)
    &=\max_{v}\left\{|x|^2_Q+|K_ix|_R^2
      +V^{ii}\left(v,Z+\begin{bmatrix}-v\\x\\-K_ix\end{bmatrix}\!\!
      \begin{bmatrix}-v\\x\\-K_ix\end{bmatrix}^\top\right)\right\}\\
    &=\max_{v}\left\{|x|^2_Q+|K_ix|_R^2+|v|_{P_{ii}}^2-\gamma^2|A_ix+B_iu-v|^2-z_i\right\}\\
    &=|x|^2_Q+|K_ix|^2_R+\left|(A_i-B_iK_i)x\right|^2_{(P_{ii}^{-1}-\gamma^{-2}I)^{-1}}-z_i\\
    &\le|x|^2_{P_{ii}}-z_i\\
    &=V^{ii}(x,Z).  
  \end{align*}
  If not, Lemma~\ref{lem:T} in the Appendix together with (\ref{eqn:Pik}) give
  \begin{align*}
    \mathcal{F}_{-K_kx}V^{ij}(x,Z)
    &=\max_{v}\left\{|x|^2_Q+|K_kx|_R^2
      +V^{ij}\left(v,Z+\begin{bmatrix}-v\\x\\-K_kx\end{bmatrix}\!\!
      \begin{bmatrix}-v\\x\\-K_kx\end{bmatrix}^\top\right)\right\}\\
    &=\max_{v}\left\{|x|^2_Q+|K_kx|_R^2+|v|_{P_{ij}}^2-\frac{1}{2}\sum_{l\in\{i,j\}}\left(\gamma^2|A_lx+B_lu-v|^2+z_l\right)\right\}\\
    &=|x|^2_Q+|K_{k}x|^2_R+\bigg|(A_i-B_iK_k+A_j-B_jK_k)x/2\bigg|^2_{(P_{ij}^{-1}-\gamma^{-2}I)^{-1}}\notag\\
    &\quad-\gamma^2|(A_i-B_iK_k-A_j+B_jK_k)x/2|^2-(z_i+z_j)/2\\
    &\le\max\{|x|^2_{P_{ik}},|x|^2_{P_{jk}}\}-(z_i+z_j)/2\\
    &=\max\{V^{ik}(x,Z)+(z_k-z_j)/2,V^{jk}(x,Z)+(z_k-z_i)/2\}\\
    &\le \max\{V^{ik}(x,Z),V^{jk}(x,Z)\}.
  \end{align*}
  Hence
  \begin{align*}
    \mathcal{F}\bar{V}(x,Z)
    &\le\mathcal{F}_{-K_kx}\bar{V}(x,Z)
    =\max_{i,j}\mathcal{F}_{-K_kx}V^{ij}(x,Z)
    \le\max_{i}V^{ik}(x,Z)\le\bar{V}(x,Z).
  \end{align*}
    so the Bellman inequality $\bar{V}\le\mathcal{F}\bar{V}$ is proved.
    Moreover, the control law $u=-K_kx$ can equivalently be written as (\ref{eqn:Kcontrol}), since (\ref{eqn:Z_t}) gives
    \begin{align*}
      \left\|\begin{bmatrix}I&A_i&B_i\end{bmatrix}^\top\right\|_{Z_t}^2
      &=\sum_{\tau=0}^{t-1}|A_ix_\tau+B_iu_\tau-x_{\tau+1}|^2.
    \end{align*}    
    Finally, the inequality $J_\mu(x_0)\le\max_{i,j}|x_0|^2_{P_{ij}}$ follows by application of  Theorem~\ref{thm:intro}. 
\end{pf}

\section{Conclusions}
Two main contributions have been given in the is paper. First, Theorem~\ref{thm:intro} shows that the minimax adaptive control problem can be stated as a standard zero sum dynamic game with a finite-dimensional state space. Second, Theorem~\ref{thm:main} gives and explicit upper bound for the value of the game and a corresponding adaptive control scheme achieving the upper bound. We believe that the approach has tremendous potential for applications and generalizations in various directions. However, we also believe that upper bound of Theorem~\ref{thm:main} may not be the last word, but rather a first step towards an exact formula for the minimal achievable gain and a seed for a rich new theory of the interplay between learning and control.

\acks{The author is grateful to Olle Kjellqvist for help with constructing the example and to Daniel Cederberg for spotting an error in the original publication of 2021.

 The author is a member of the excellence center ELLIIT. Financial support was obtained from the Swedish Research Council, the European Research Council and the Swedish Foundation for Strategic Research. The work was also partially supported by the Wallenberg AI, Autonomous Systems and Software Program (WASP) funded by the Knut and Alice Wallenberg Foundation.
}

\section{Appendix: Quadratic optimization.}
\begin{lem}
  For $P\in\realR^{n\times n}$, $y,v\in\realR^n$ and $\gamma\in\realR$, it holds that
  \begin{align*}
    \max_{v}\left\{|v|^2_P-\gamma^2|y-v|^2\right\}
    &=|y|^2_{(P^{-1}-\gamma^{-2}I)^{-1}}.
  \end{align*}
\label{lem:P}
\end{lem}
\begin{pf}
  \begin{align*}
    \max_{v}\left\{|v|^2_P-\gamma^2|y-v|^2\right\}
    &=\max_{v}\left\{|v|^2_{(P-\gamma^2I)}-\gamma^2|y|^2+2\gamma^2y^\top v\right\}\\
    &=-\gamma^2|y|^2-\gamma^4y^\top(P-\gamma^2I)^{-1}y\\
    &=\gamma^2y^\top(P-\gamma^2I)^{-1}\left[\gamma^2I-P-\gamma^2I\right]y\\
    &=y^\top(P^{-1}-\gamma^{-2}I)^{-1}y
  \end{align*}
\end{pf}

\begin{lem}
  For $T\in\realR^{n\times n}$, $y_1,y_2,v\in\realR^n$ and $\gamma\in\realR$, it holds that
  \begin{align*}
    &\max_{v}\left\{|v|^2_T-{\gamma^2}|y_1-v|^2/2-\gamma^2|y_2-v|^2/2\right\}\\
    &=|(y_1+y_2)/2|^2_{(T^{-1}-\gamma^{-2}I)^{-1}}-\gamma^2|(y_1-y_2)/2|^2
  \end{align*}
\label{lem:T}
\end{lem}
\begin{pf}
Application of Lemma~\ref{lem:P} gives
  \begin{align*}
    &\max_{v}\left\{|v|^2_T-{\gamma^2}|y_1-v|^2/2-\gamma^2|y_2-v|^2/2\right\}\\
    &=\max_{v}\left\{|v|^2_T-{\gamma^2}|(y_1+y_2)/2-v|^2
    +\gamma^2|(y_1+y_2)/2|^2-\gamma^2(|y_1|^2+|y_2|^2)/2\right\}\\
    &=|(y_1+y_2)/2|^2_{(T^{-1}-\gamma^{-2}I)^{-1}}
    +\gamma^2|(y_1+y_2)/2|^2-\gamma^2(|y_1|^2+|y_2|^2)/2\\
    &=|(y_1+y_2)/2|^2_{(T^{-1}-\gamma^{-2}I)^{-1}}-\gamma^2|(y_1-y_2)/2|^2.
  \end{align*}
\end{pf}


\end{document}